\documentclass[12pt]{amsart}
\numberwithin{equation}{section}

\usepackage{calc}  
\usepackage{amsmath}  
\usepackage{amsthm}
\usepackage{amsfonts}
\usepackage{amsopn}
\usepackage{amssymb}
\usepackage{amstext}
\usepackage{amscd}
\usepackage{graphicx}
\usepackage[all]{xy}

\newtheorem{theorem}{Theorem}[section]
\newtheorem{proposition}[theorem]{Proposition}
\newtheorem{corollary}[theorem]{Corollary}
\newtheorem{lemma}[theorem]{Lemma}

\theoremstyle{definition}

\newtheorem{example}[theorem]{Example}
\newtheorem{definition}[theorem]{Definition}
\newtheorem{remark}[theorem]{Remark}
\newtheorem{acknow}{Acknowledgments}

\def\Z{\mathbb{Z}}

\def\R{\mathbb{R}}
\def\C{\mathbb{C}}

\def\CP{\mathbb{C} \mathbb{P}}

\def\<{\left\langle}
\def\>{\right\rangle}
\def\({\left(}
\def\){\right)}

\def\cC{\mathcal{C}}
\def\cD{\mathcal{D}}

\def\cF{\mathcal{F}}
\def\cG{\mathcal{G}}
\def\cH{\mathcal{H}}

\def\cL{\mathcal{L}}

\def\cM{\mathcal{M}}
\def\cN{\mathcal{N}}
\def\cO{\mathcal{O}}
\def\cP{\mathcal{P}}

\def\cY{\mathcal{Y}}

\def\fs{\mathfrak{s}}

\def\Conf{\operatorname{Conf}}

\def\Ker{\operatorname{Ker}}
\def\Im{\operatorname{Im}}

\def\Ind{\operatorname{Ind}}

\begin{document}

\title[Spin Structures on Seiberg-Witten Moduli Spaces]
{Spin Structures on \\ Seiberg-Witten Moduli Spaces}

\author[H. Sasahira]{Hirofumi Sasahira$^*$}

\renewcommand{\thefootnote}{\fnsymbol{footnote}}
\footnote[0]{2000\textit{ Mathematics Subject Classification}.
Primary 57R57; Secondary 53C25.}

\thanks{
$^*$Partially supported by the 21th century COE program
at Graduate School of Mathematical Sciences,
the University of Tokyo.
}

\address{
Graduate School of Mathematical Sciences,
University of Tokyo,\endgraf 
3-8-1 Komaba Meguro-ku, Tokyo 153-8941, Japan
}
\email{sasahira@ms.u-tokyo.ac.jp}

\date{}
\maketitle

\begin{abstract}
Let $M$ be an oriented closed $4$-manifold with a spin$^c$ structure $\cL$.
In this paper we prove that  
under a suitable condition for $(M,\cL)$
the Seiberg-Witten moduli space has a canonical spin structure
and its spin bordism class is an invariant of $M$.
We show that the invariant of $M=\#_{j=1}^l M_j$ is non-trivial for some spin$^c$ structure when $l$ is $2$ or $3$ and each $M_j$ is a $K3$ surface or 
a product of two oriented closed surfaces of odd genus. 
As a corollary, we obtain the adjunction inequality for the $4$-manifold.
Moreover we calculate the Yamabe invariant of $M \# N_1$ for some negative definite $4$-manifold $N_1$.
We also show that
$M \# N_2$ does not admit an Einstein metric for some negative definite $4$-manifold $N_2$.

\end{abstract}

\section{Introduction}
Since E. Witten introduced the Seiberg-Witten equations (\cite{W}),
the moduli space of solutions to the equations
has been applied to $4$-dimensional topology.
M. Furuta used the Seiberg-Witten equations themselves, rather than the moduli space, to obtain the $10/8$ theorem (\cite{F}).
Roughly speaking, the Seiberg-Witten moduli space is the zero locus of the map defining the equations,
which we call the Seiberg-Witten map, 
between two Hilbert bundles over the Jacobian torus.
Furuta used finite dimensional approximation of the Seiberg-Witten map
to prove the $10/8$ theorem.
Moreover using finite dimensional approximation of the Seiberg-Witten map,
S. Bauer and Furuta refined the Seiberg-Witten invariants (\cite{BF}).
The refined invariant is more powerful than the usual Seiberg-Witten invariant.
There are $4$-manifolds for which
the usual Seiberg-Witten invariants vanish but
the Bauer-Furuta invariants do not
(\cite{B,FKM}).
It is, however, hard in general to detect the Bauer-Furuta invariants. 

To detect the Bauer-Furuta invariants explicitly,
we define new invariants for $4$-manifolds.
This invariant is weaker than the Bauer-Furuta invariant, but easier to treat, in particular when the first Betti number of the $4$-manifold is positive.
An outline of the definition of the invariant is as follows.

Let $(M,g)$ be an oriented, closed $4$-dimensional Riemannian manifold 
with $b^+(M)>1$,
and $\cL$ a spin$^c$ structure on $M$.
We write $\Ind (D)$ 
for the index bundle of the Dirac operators parameterized by 
$T=H^1(M;\R)/H^1(M;\Z)$
(see \S \ref{sufficient condition}) .
If $c_1(\Ind (D)) \equiv 0 \mod 2$, 
then the Seiberg-Witten moduli space allows a spin structure,
and a choice of square root of the determinant line bundle $\det \Ind (D)$ 
determines a spin structure of the moduli space.
The spin bordism class of the moduli space is an invariant of $M$ 
which depends only on $\cL$ and the choice of square root of $\Ind(D)$.

We calculate the invariant for $M=\#_{j=1}^l M_j$ when $M_j$ is a $K3$ surface or 
a product of two oriented closed surfaces of odd genus,
and $l$ is $2$ or $3$.
We take a spin$^c$ structure on $M$ of the form $\cL=\#_{j=1}^l\cL_j$,
where $\cL_j$ is a spin$^c$ structure on $M_j$ 
induced by a complex structure.
We show that 
in this case $c_1(\Ind(D)) \equiv 0 \mod 2$
and our invariant is non-trivial.
As an application, we obtain the adjunction inequality for such $M$,
i.e., if an oriented closed surface $\Sigma$ of positive genus  is  embedded in $M$ satisfying that its self-intersection number $\Sigma \cdot \Sigma$ is nonnegative, then we have
\[
\Sigma \cdot \Sigma
\leq
\< c_1(\det \cL),\Sigma \>+2g(\Sigma)-2.
\]
Here $\det \cL$ is the determinant complex line bundle of $\cL$,
and $g(\Sigma)$ is the genus of $\Sigma$.

As another application, following Ishida and LeBrun's argument in \cite{IL},
we compute the Yamabe invariant of $M \# N_1$ when $N_1$ is an oriented, closed, negative definite $4$-manifold admitting a Riemannian metric with scalar curvature nonnegative at each point.
We also show that
if $N_2$ is an oriented, closed, negative definite $4$-manifold satisfying
\[
4l-(2\chi(N_2)+3\tau (N_2)) \geq 
\frac{1}{3}\sum_{j=1}^l c_1(M_j)^2,
\]
then $M \# N_2$ does not admit an Einstein metric,
where $\tau(N_2)$ and $\chi(N_2)$ are the signature and the Euler number of $N_2$ respectively.

\begin{acknow}
This paper is part of the author's master thesis.
The author would like to thank my advisor Mikio Furuta 
for his suggestions and warm encouragement.
The author also thanks Masashi Ishida for useful information
about Einstein metrics and Yamabe invariants.
\end{acknow}

\section{Finite dimensional approximations of the Seiberg-Witten map}

In this section, we review the definition of the Seiberg-Witten map and 
its finite dimensional approximation.
We refer the readers to \cite{BF} for details.

\subsection{The Seiberg-Witten map}

Let $M$ be an oriented, closed, connected $4$-manifold 
and let $g$ be a Riemannian metric on $M$.
Assume that $b^+(M)>1$.
Choose a spin$^c$ structure $\cL$ on $M$.
We write $S^{\pm}(\cL)$ for the positive and negative spinor bundles, and $\det \cL$ for the determinant line bundle 
associated with $\cL$.

Let $k$ be an integer larger than or equal to $4$
and set 
$\hat{\cG}=\{ \gamma \in L^2_{k+1}(M,U(1))|\gamma(x_0)=1 \}$
for a fixed base point $x_0 \in M$.
Fix a connection $A_0$ on $\det \cL$,
and define $T:=(A_0+i\Ker d)/\hat{\cG}$,
where $d:L^2_k(T^*M) \rightarrow L^2_{k-1}(\Lambda^2 T^*M)$ is 
the exterior derivative.
The action of $\gamma \in \hat{\cG}$ on $A \in (A_0+i\Ker d)$ is defined by
\begin{equation} \label{gauge}
\gamma A :=A+2\gamma^{-1} d\gamma.
\end{equation}
Put
\begin{gather*}
\tilde{\cC}(\cL):=L^2_k(S^+(\cL) \oplus T^*M), \\
\tilde{\cY}(\cL):=
L^2_{k-1}(S^-(\cL) \oplus \Lambda^+T^*M)
\oplus \cH_g^1(M) \oplus (L^2_{k-1}(M)/\R),
\end{gather*}
where $\R$ represents the space of constant functions on $M$
and $\cH_g^1(M)$ is the space of harmonic $1$-forms on $M$ with respect to $g$.
Let $\cC(\cL) \rightarrow T$ and $\cY(\cL) \rightarrow T$ be
Hilbert bundles on $T$ defined by
\begin{gather*}
\cC(\cL):=(A_0+i\Ker d) \times_{\hat{\cG}} \tilde{\cC}(\cL),  \\
\cY(\cL):=(A_0+i\Ker d) \times_{\hat{\cG}} \tilde{\cY}(\cL).
\end{gather*}
The action of $\hat{\cG}$ 
on $(A_0+i\Ker d)$ is given by $(\ref{gauge})$.
We define actions of $\hat{\cG}$
on $L^2_k(S^+(\cL))$ and 
on $L^2_{k-1}(S^-(\cL))$ by fiber-wise scalar products.
We define actions of $\hat{\cG}$ on the other terms to be trivial.
We define $U(1)$-actions on $\cC(\cL)$ and $\cY(\cL)$ by scalar products 
on $L^2_k(S^+(\cL))$ and $L^2_{k-1}(S^-(\cL))$ and set
\[
\cP:=\left\{ (g,\eta) \in Riem(M) \times L^2_k(\Lambda^2 T^*M) 
\left| 
[\eta]_g^+ \not= [F_{A_0}]_g^+
\right. \right\},
\]
where $Riem(M)$ is the space of Riemannian metrics on $M$,
and $[\eta]_g^+$ and $[F_{A_0}]_g^+$ are 
$\cH_g^+(M)$ parts of $\eta$ and $F_{A_0}$ respectively.
For $(g,\eta) \in \cP$, we define the Seiberg-Witten map by
\[ \begin{array}{rccc}
SW_{g,\eta}: &\cC(\cL) &\longrightarrow& \cY(\cL) \\
& (A,\phi,a) &\longmapsto& 
(A,D_{A+ia}\phi,F^+_{A+ia} -q(\phi)-\eta^+ ,p(a),d^*a),
\end{array} \]
where $q(\phi)$ is a quadratic form of $\phi$ and
$p:L^2_k(T^*M) \rightarrow \cH_g^1(M)$ is the $L^2$-projection.
The moduli space $\cM_M(\cL,g,\eta)$ 
of solutions to the Seiberg-Witten equations perturbed by $(g,\eta)$
is identified with $SW_{g,\eta}^{-1}(0)/U(1)$ naturally.

The following fact is well known.
\begin{theorem} [\cite{KM}] \label{compactness}
For generic $(g,\eta) \in \cP$,
$\cM_M(\cL,g,\eta)$ is a compact smooth manifold
and an orientation on $\cH_g^1(M;\R) \oplus \cH_g^+(M;\R)$
determines an orientation on $\cM_M(\cL,g,\eta)$.
\end{theorem}

\subsection{Finite dimensional approximation}
We explain finite dimensional approximations 
of the Seiberg-Witten map briefly.

Let $\cD:\cC(\cL) \rightarrow \cY(\cL)$
be the linear part of the SW map:
\[ \begin{array}{rccc}
\cD: &\cC(\cL) &\longrightarrow& \cY(\cL) \\
& (A,\phi,a) &\longmapsto& (A, D_A \phi, d^+a, p(a), d^*a).
\end{array} \]
By Kuiper's theorem \cite{Ku},
we have a global trivialization of $\cY(\cL)$
\[
\cY(\cL) \cong T \times H,
\]
where $H$ is a Hilbert space.
We fix a trivialization of $\cY(\cL)$.
Since $\cY(\cL)$ has the complex part and the real part,
$H$ decomposes into the direct sum $H_{\C} \oplus H_{\R}$
of a complex Hilbert space $H_{\C}$ and
a real Hilbert space $H_{\R}$.

For a finite dimensional subspace $W \subset H$,
let $p_{W}:\cY(\cL)=T \times H \rightarrow W$ 
be the projection.
We denote $\cD^{-1}(T \times W)$ by $\cF(W)$.
Then we define $f_W:\cF(W) \rightarrow W$ by
\[
f_W=p_W \circ SW|_{\cF(W)}:\cF(W) \longrightarrow W.
\]

\begin{theorem} [\cite{BF}] \label{approximation}
Let $W^+$ and $\cF(W)^+$ be 
the one-point compactifications of $W$ and $\cF(W)$.
Then $f_W:\cF(W) \rightarrow W$ induces a $U(1)$-equivariant map
$f_{W}^+:\cF(W)^+ \rightarrow W^+$,
and there is a finite dimensional subspace $W \subset H$ 
such that the following conditions are satisfied. \\

$(1)$
$\Im \cD + (T \times W)=\cY(\cL)$. \\

$(2)$
For all finite dimensional subspace $W' \subset H$
such that $W \subset W'$,
the diagram
\[
\xymatrix{
\quad \quad \cF(W')^+ \quad  \quad \ar[rrr]^{\quad f_{W'}^+} \ar@{=}[d] & & & 
\quad (W')^+ \quad \ar@{=}[d] \\
(\cF(W) \oplus \cF(U))^+ \ar[rrr]_{\ \ \ (f_W \oplus p_U \cD|_{\cF(U)})^+} & & & 
(W \oplus U)^+
}
\]
is $U(1)$-equivariant homotopy commutative as pointed maps,
where $U$ is the orthogonal complement of $W$ in $W'$.

\end{theorem}

\begin{definition}
When $W \subset H$ satisfies $(1)$ and $(2)$,
we call $f_W:\cF(W) \rightarrow W$ 
a finite dimensional approximation of the Seiberg-Witten map.
\end{definition}

\section{Spin structures on moduli spaces}
In \S \ref{sufficient condition},
by using finite dimensional approximation of the Seiberg-Witten map,
we show a sufficient condition for the moduli space to be a spin manifold.
In \S \ref{invariant},
we will prove that
the spin bordism class of the spin structure on the moduli space
is an invariant of $M$.
In \S \ref{example}, we give some applications of this invariant.

\subsection{A sufficient condition for moduli space to have a spin structure.} 
\label{sufficient condition}

Let $f=f_W:V=\cF(W) \rightarrow W$ 
be a finite dimensional approximation of the Seiberg-Witten map.
Note that $V$ has a natural decomposition $V=V_{\C} \oplus V_{\R}$
into the direct sum of 
a complex vector bundle and a real vector bundle.
Similarly decompose $W$ as $W=W_{\C} \oplus W_{\R}$.

If we take a generic $(g,\eta) \in \cP$ as in Theorem \ref{compactness},
$\cM_M(\cL,g,\eta)$ does not include reducible monopoles,
hence $f^{-1}(0)$ lies in $V_{irr}:=(V_{\C} \backslash \{ 0 \}) \times_T V_{\R}$.
Put $\bar{V}:=V_{irr}/U(1)$ 
and $\cM:=f^{-1}(0)/U(1)$.
We define a vector bundle $\bar{E} \rightarrow \bar{V}$ by
$\bar{E}:=V_{irr} \times_{U(1)} W=\bar{E}_{\C} \oplus \bar{E}_{\R}$,
where $\bar{E}_{\C}=V_{irr} \times_{U(1)} W_{\C}$, $\bar{E}_{\R}=V_{irr} \times W_{\R}$.
Since $f$ is $U(1)$-equivariant,
$f$ induces a section $s:\bar{V} \rightarrow \bar{E}$.
Then $\cM$ is the zero locus of $s$.
If necessary, we perturb $s$ on a compact subset in $\bar{V}$ so that
$s$ is transverse to the zero section of $\bar{E}$
and $\cM$ is a compact smooth submanifold of $\bar{V}$.

We can orient $\cM$ by using an orientation on $\cH_g^1(X) \oplus \cH_g^+(X)$ in the following way.
The real part $\cD_{\R}$ of $\cD$ is independent of $A \in T$ and the cokernel is naturally identified with $\cH_g^+(X)$. So $W_{\R}$ has the form $\cH_g^+(X) \oplus W_{\R}'$ and $\cD_{\R}$ induces the isomorphism between each fiber of $V_{\R}$ and $W_{\R}'$.
(Hence $V_{\R}$ is a trivial vector bundle on $T$.)
If we choose orientations on $W_{\R}'$ and $\cH_g^+(X)$,
we get an orientation on $\bar{E}_{\R}$ and 
orientations on $V_{\R}$ and $\cH_g^1(X)$ compatible with $\cD_{\R}$ and $\cO$.
$T$ is naturally identified with $H^1(X;\R)/H^1(X;\Z)$, so
the tangent bundle $T(T)$ of $T$ has a natural trivialization $T(T) \cong T \times H^1(X;\R) \cong T \times \cH_g^1(X)$. The orientation on $\cH_g^1(X)$ induces an orientation on $T(T)$. 
These orientations induce an orientation on the tangent bundle $T\bar{V}$ by Lemma \ref{tangent bar V} below.
The derivative of $s$ induces an isomorphism between $\bar{E}|_{\cM}$ and the normal bundle $\cN$ of $\cM$ in $\bar{V}$.
The orientation on $\bar{E}$ induces an orientation on $\cN$ through this isomorphism, and we have an orientation on $\cM$ compatible with the decomposition
$T\bar{V}|_{\cM} = T\cM \oplus \cN$.
(It is easy to check that this orientation on $\cM$ is independent of the choices of the orientations on $W_{\R}'$ and $\cH_g^+(X)$.)
So we have the following.

\begin{lemma}
A choice of orientation on $\cH_g^1(X) \oplus \cH_g^+(X)$ induces an orientation on $\cM$.
\end{lemma}

When $T\bar{V}$ and $\bar{E}$ have spin structures, we can equip $\cM$ with  a spin structure as in the case of orientation.
The spin structure on $\bar{E}$ induces a spin structure on $\cN$ through the derivative of $s$.
Since $T\bar{V}|_{\cM}$ is the direct sum of 
$T\cM$ and $\cN$,
spin structures on $T\bar{V}$ and $\cN$ induce a spin structure on $\cM$,
from the next well-known lemma.

\begin{lemma} \label{spin sum}
Let $X$ be a smooth manifold,
$F_1$ and $F_2$ be oriented vector bundles on $X$.
If $F_1$ and $F_2$ have spin structures,
then spin structures on $F_1$ and $F_2$ determine a spin structure on 
$F_1 \oplus F_2$.
If $F_1$ and $F_1 \oplus F_2$ have spin structures,
then spin structures on $F_1$ and $F_1 \oplus F_2$ determine 
a spin structure on $F_2$ naturally.
\end{lemma}

Therefore we have shown the following.

\begin{lemma} \label{key lemma}
Let $f:V \rightarrow W$ be a finite dimensional approximation
of the Seiberg-Witten map.
Assume that $T\bar{V}$ and $\bar{E}$ have a spin structure.
Choose spin structures $\fs_{\bar{V}}$ and $\fs_{\bar{E}}$
on $T\bar{V}$ and $\bar{E}$.
Then $\fs_{\bar{V}}, \fs_{\bar{E}}$ and $f$ induce 
a spin structure on $\cM=f^{-1}(0)/U(1)$.
\end{lemma}

We calculate $w_2(T\bar{V})$ and $w_2(\bar{E})$ 
to know when $T\bar{V}$ and $\bar{E}$ have spin structures.

Let $a \in \Z$ be the index of the Dirac operator,
let $\Ind D \in K(T)$ be the index bundle of the Dirac operators 
$\{ D_A \}_{A \in T}$ parameterized by $T$.
Then we have $\Ind D = [V_{\C}] - [\underline{\C}^m] \in K(T)$, $V_{\R}=\underline{\R}^n$, $W_{\C}=\C^m$, $W_{\R}=\cH_g^+(X) \oplus W_{\R}'$, $\dim W_{\R}'=n$
for some $m,n \in \Z_{\geq 0}$.

\begin{lemma} \label{tangent bar V}
Let $\bar{\pi} :\bar{V} \rightarrow T$ be the projection
and define a complex line bundle $H \rightarrow \bar{V}$ by
$H:=V_{irr} \times_{U(1)} \C$.
Then there is a natural isomorphism 
\[
T\bar{V} \oplus \underline{\R} 
\cong
\bar{\pi}^*T(T) \oplus (\bar{\pi}^*V_{\C} \otimes_{\C} H) 
\oplus \bar{\pi}^*V_{\R}.
\]
\end{lemma}

\begin{proof}
Let $\pi_{irr}:V_{irr} \rightarrow T$ and 
$p:V_{irr} \rightarrow \bar{V}=V_{irr}/U(1)$ be the projections.
Note that we have a $U(1)$-equivariant isomorphism
\[
p^*(T\bar{V}) \oplus \underline{\R}
\cong
TV_{irr}=\pi_{irr}^*(T(T) \oplus V).
\]
where $\underline{\R}$ stands for the $U(1)$-orbit direction.
Then by dividing by the $U(1)$-actions,
we obtain the required isomorphism.
\end{proof}

By Lemma $\ref{tangent bar V}$ and the triviality of $V_{\R}$, 
we have
$w_2(T\bar{V}) \equiv 
\bar{\pi}^* c_1(V_{\C})+(m+a)c_1(H)
\mod 2$.
By $(1)$ in Theorem \ref{approximation},
$c_1(V_{\C})$ is equal to $c_1(\Ind(D))$,
thus we have
\begin{equation} \label{barV D H}
w_2(T\bar{V}) 
\equiv 
\bar{\pi}^*c_1(\Ind(D))+(m+a)c_1(H)
\mod 2 .
\end{equation}
T-J. Li and A. Liu calculated
$c_1(\Ind (D))$ in \cite{LL} as follows.

Let $\{ \alpha_j \}_{j=1}^{b_1}$ be generators of $H^1(M;\Z)$.
Then we obtain a natural identification,
\[
T \cong H^1(M;\R)/H^1(M;\Z)
\cong \R^{b_1}/\Z^{b_1}=T^{b_1}.
\]
Let $\Psi$ be a map $M \rightarrow T^{b_1} \cong T$ given by
\[
x \longmapsto
\left( \int_{x_0}^x \alpha_1, \cdots, \int_{x_0}^x \alpha_{b_1} \right).
\]
This map is well defined by the Stokes theorem
and induces the isomorphism
$\Psi^*:H^1(T;\Z) \cong H^1(M;\Z)$.
Set $\beta_j=(\Psi^*)^{-1}(\alpha_j) \in H^1(T;\Z)$.

\begin{proposition} [\cite{LL}] \label{ind D}
Let $\Ind D \in K(T)$
 be the index bundle of the Dirac operators $\{ D_A \}_{A \in T}$
parameterized by $T$.
Then the first Chern class $c_1(\Ind (D))$ of $\Ind(D)$ is given by
\[
c_1(\Ind (D))=
\frac{1}{2}\sum_{i<j} 
\< c_1(\det \cL) \alpha_i \alpha_j,[M] \> \beta_i \beta_j
\in H^2(T;\Z).
\]
\end{proposition}

By the equation $(\ref{barV D H})$ and Proposition $\ref{ind D}$,
we have the following.

\begin{lemma} \label{w2V}
The second Stiefel-Whiteny class of $T\bar{V}$ is given by
\[
w_2(T\bar{V}) \equiv
\sum_{i<j} c_{ij} \bar{\pi}^* \beta_i \beta_j+(m+a)c_1(H) \mod 2,
\]
where $\displaystyle 
c_{ij}:=\frac{1}{2}\< c_1(\det \cL)\alpha_i \alpha_j,[M] \>$.
\end{lemma}

On the other hand,
by the definitions of $\bar{E}$ and $H$,
we have $\bar{E}=mH \oplus \underline{\R}^{n+b}$.
Hence we obtain the following.

\begin{lemma} \label{w2E}
The second Stiefel-Whiteny class of $\bar{E}$ is given by
\[
w_2(\bar{E}) \equiv mc_1(H) \mod 2.
\]
\end{lemma}

By Lemma \ref{key lemma}, Lemma \ref{w2V} and Lemma \ref{w2E},
we have the following.

\begin{proposition} \label{the conditions}
Let
$f:V \rightarrow W$ be a finite dimensional approximation
of the Seiberg-Witten map such that
$m=\dim_{\C} W_{\C}$ is even.
Then $T\bar{V}$ and $\bar{E}$ have a spin structure
if the pair $(M,\cL)$ satisfies the following conditions.
\[
(*) \left\{
\begin{array}{lrl}
(*)_1 & a \equiv 0      &\mod 2 \\
(*)_2 & c_{ij} \equiv 0 &\mod 2~ (\forall i,j).
\end{array} 
\right.
\]
Moreover if we choose spin structures $\fs_{\bar{V}}$ and $\fs_{\bar{E}}$
of $T\bar{V}$ and $\bar{E}$,
then $\fs_{\bar{V}}$, $\fs_{\bar{E}}$ and $f$  
equip $\cM$ with a spin structure.
\end{proposition}

\subsection{Invariants for $4$-manifolds defined by spin structures on $\cM$}
\label{invariant}

An orientation on $\cH_g^1(M) \oplus \cH_g^+(M)$ determines
an orientation on $\cM$ (\S \ref{sufficient condition}). 
We will show that
when the condition $(*)$ is satisfied,
a certain datum in addition to the orientation on $\cH_g^1(M) \oplus \cH_g^+(M)$ determines a canonical spin structure on $\cM$.
The datum is actually a square root of $\det \Ind(D)$.
To explain it,
we need the following lemma.

\begin{lemma} \label{root}
Let $X$ be a smooth manifold and
$F \rightarrow X$ be a complex bundle with $c_1(F) \equiv 0 \mod 2$.
A choice of complex line bundle 
$L \rightarrow X$ which satisfies $L^{\otimes 2}=\det F$
naturally determines a spin structure on $F$.
\end{lemma}

\begin{proof}
The $2$-fold cover of $U(n)$ is given by
\[
\{ (A,t) \in U(n) \times S^1|\det A=t^2 \},
\]
which is naturally regarded as a subgroup of $Spin(2n)$.
Take an open covering $\{ U_j \}_j$ of $X$ such that
$F$ and $L$ have trivializations on each $U_j$.
Fix hermitian metrics on $F$ and $L$
compatible with the identification $L^{\otimes 2}=\det F$.
We denote transition functions on $U_i \cap U_j$ of $F$ and $L$ by
$g_{ij}:U_i \cap U_j \rightarrow U(n)$ and
$z_{ij}:U_i \cap U_j \rightarrow S^1$.
Then $\det g_{ij}=z_{ij}^2$,
since $\det F=L^{\otimes 2}$.
Put $\tilde{g}_{ij}=(g_{ij},z_{ij}):U_i \cap U_j \rightarrow Spin(2n)$,
then $\{ \tilde{g}_{ij} \}_{ij}$ satisfies the cocycle condition and
determines a spin structure of $F$.
\end{proof}

When the condition $(*)_2$ is satisfied,
then $c_1(\Ind (D)) \equiv 0 \mod 2$.
So we can take a complex line bundle 
$L \rightarrow T$ such that 
$L^{\otimes 2}=\det \Ind (D)$.

\begin{proposition}
Assume that the pair $(M,\cL)$ satisfies the conditions $(*)$.
Let $f:V \rightarrow W$ be a finite dimensional approximation of 
the Seiberg-Witten map such that $m=\dim_{\C}W_{\C}$ is even.
Then the finite dimensional approximation $f$, 
an orientation $\cO$ of $\cH_g^1(M) \oplus \cH_g^+(M)$ and 
a complex line bundle
$L \rightarrow T$ which satisfies $L^{\otimes 2}= \det \Ind (D)$
determine a canonical spin structure on $\cM$.
\end{proposition}

\begin{proof}
Suppose that the pair $(M,\cL)$ satisfies the condition $(*)$.
By Lemma \ref{key lemma},
spin structures on $T\bar{V}$, $\bar{E}$ 
and a finite dimensional approximation $f$
induce a canonical spin structure on $\cM$.
So it is sufficient to show that
$\cO$ and $L$ induce spin structures on $T\bar{V}$ and $\bar{E}$.
By Lemma \ref{tangent bar V},
we have only to show that
the choices of $\cO$ and $L$ induce spin structures on
$\bar{\pi}^* V_{\C} \otimes H$,  $V_{\R}$, $T(T)$ and $\bar{E}$.

Since $m$ is even and condition $(*)_1$ is satisfied,
$\bar{\pi}^*L \otimes H^{\otimes \frac{m+a}{2}}$ is a square root of
$\det (\bar{\pi}^*V_{\C} \otimes H)=
(\bar{\pi}^* \det V_{\C}) \otimes H^{\otimes (m+a)}$.
So by Lemma $\ref{root}$,
we have a spin structure on $\bar{\pi}^* V_{\C} \otimes H$.

Recall that $W_{\R}$ is the direct sum $\cH_g^+(X) \oplus W_{\R}'$.
We fix orientations on $\cH_g^+(X)$ and $W_{\R}'$, then we have orientations on $V_{\R}$ and $\cH_g^1(X)$ compatible with $\cD_{\R}$ and $\cO$.
(See \S \ref{sufficient condition}.)
Since the real part $\cD_{\R}$ of $\cD$ is independent of $A \in T$, $V_{\R}$ has a natural trivialization compatible with the orientation.
This trivialization equips $V_{\R}$ with a spin structure.
The tangent bundle $T(T)$ of $T$ has a natural trivialization $T(T) = T \times \cH_g^1(M)$ and the orientation $\cH_g^1(X)$ orients $T(T)$.
So we have a spin structure on $T(T)$ compatible with this trivialization.

Lastly we consider $\bar{E}$.
Let $\bar{E}_{\C}$ be the complex part of $\bar{E}$,
i.e. $\bar{E}_{\C}=V_{irr} \times_{U(1)} \C^m$.
Since $\det \bar{E}_{\C}=H^{\otimes m}$,
$H^{\otimes \frac{m}{2}}$ is a square root of $\det \bar{E}_{\C}$.
So by Lemma $\ref{root}$,
a spin structure of $\bar{E}_{\C}$ is determined.
Let $\bar{E}_{\R}$ be the real part of $\bar{E}$.
Then $\bar{E}_{\R}=V_{irr} \times W_{\R} = V_{irr} \times (\cH_g^+(X) \oplus W_{\R}')$.
Hence $\bar{E}_{\R}$ has a natural spin structure
induced by the trivialization.

We have seen that $f$, $\cO$ and $L$ determine a spin structure on $\cM$ if we choose orientations on $\cH_g^+(X)$ and $W_{\R}'$.
It is easy to see that this spin structure is independent of the choices of orientations on $\cH_g^+(X)$ and $W_{\R}'$.
\end{proof}

Let $\pi: \cM \rightarrow T$ be
the restriction of the projection $\bar{V} \rightarrow T$ to $\cM$.
We show that the class
$(\cM, \pi) \in \Omega^{spin}_d(T)$ induced by $f,\cO,L$ 
is an invariant of $M$.
Here $d$ is the dimension of $\cM$.

\begin{theorem}
Assume that 
the pair $(M,\cL)$ satisfies the condition $(*)$.
The class $( \cM,\pi) \in \Omega^{spin}_d(T)$
which is induced by $f,\cO,L$
is independent of the perturbation $(g,\eta) \in \cP$ and 
the finite dimensional approximation $f$.
\end{theorem}

\begin{proof}
Fix $(g,\eta) \in \cP$,
and take different finite dimensional approximations
$f_i:V_i \rightarrow W_i,(i=0,1)$
of the Seiberg-Witten map $SW_{g,\eta}$.
Denote $f_i^{-1}(0)/U(1)$ by $\cM_i$
and let $\pi_i$ be the restriction of the projections
$\bar{V}_i \rightarrow T$ to $\cM_i$.
By considering a larger finite dimensional approximation
$f:V \rightarrow W$
with $V_i \subset V$ and $W_i \subset W$,
we can assume that $V_0 \subset V_1,W_0 \subset W_1$
without loss of generality.

Let $V_1=V_0 \oplus V'$ and $W_1=W_0 \oplus W'$,
then $\cD|_{V'}$ induces an isomorphism
$V' \cong T \times W'$.
By Theorem $\ref{approximation}$, the maps
\[
(f_1)^+, \ (f_0 \oplus p_{W'} \circ \cD|_{V'})^+:
V_1^+=(V_0 \oplus V')^+ \rightarrow W_1^+=(W_0 \oplus W')^+
\]
are $U(1)$-equivariantly homotopic each other as pointed maps.
It is clear that
the spin structure on $\cM_0$ induced by $f_0 \oplus p_{W'} \circ \cD|_{V'}$
is equal to one induced by $f_0$.
Let $h:[0,1] \times V_1^+ \rightarrow W_1^+$ be a homotopy
from $(f_0 \oplus \cD)^+$ to $f_1^+$ and set
$\widetilde{\cM}:=h^{-1}(0)/U(1)$.
Let $\tilde{\pi}$ be the restriction of the projection
$\bar{V}_1 \times [0,1] \rightarrow T$ to $\widetilde{\cM}$.
By using a parallel argument to 
introduce spin structures on $\cM_0$ and $\cM_1$,
we can equip $\widetilde{\cM}$ with a spin structure 
by using $h,\cO$ and $L$.
Then $(\widetilde{\cM},\tilde{\pi})$ 
is a spin bordism between $(\cM_0,\pi_0)$ and $(\cM_1,\pi_1)$.
This implies that 
when $(g,\eta) \in \cP$ is fixed,
the class $(\cM,\pi) \in \Omega^{spin}_d(T)$ is independent of 
a choice of $f$.

Next choose two elements $(g_0,\eta_0),(g_1,\eta_1) \in \cP$.
By the assumption $b^+(M)>1$,
$\cP$ is path connected,
and there is a path $(g(t),\eta(t))_{0 \leq t \leq 1}$ in $\cP$ satisfying
$(g(i),\eta(i))=(g_i,\eta_i),(i=0,1)$.
We define parameterized Seiberg-Witten map 
\[
\widetilde{SW}:[0,1] \times \cC(\cL) \rightarrow [0,1] \times \cY(\cL)
\]
in the obvious way.
Let $\tilde{f}:\widetilde{V} \rightarrow \widetilde{W}$ be
a finite dimensional approximation of $\widetilde{SW}$.
We can endow $\widetilde{\cM}=\widetilde{f}^{-1}(0)/U(1)$ with 
a spin structure in the same way as in the case of $\cM$.
Denote 
$\widetilde{V}|_{\{ i \} \times T}$ and $\widetilde{W}|_{\{ i \} \times T}$ by
$V_i$ and $W_i$ for $i=0,1$.
Since
$f_i:=\tilde{f}|_{V_i}:V_i \rightarrow W_i$ 
is a finite dimensional approximation of $SW_{g_i,\eta_i}$,
$(\widetilde{\cM},\tilde{\pi})$ is a bordism between
$(\cM_0,\pi_0)$ and $(\cM_1,\pi_1)$.
It is showed that
the class $(\cM,\pi) \in \Omega^{spin}_d(T)$ is independent of a choice of 
$(g,\eta) \in \cP$.
\end{proof}

\begin{definition}
We write $\sigma_M(\cL,\cO,L)$ for
the class in $\Omega^{spin}_d(T)$ represented by the spin structure 
on $\cM$ induced by $f,\cO,L$
and the restriction $\pi$ of 
the projection $\bar{V} \rightarrow T$ to $\cM$.
Here $d$ is the dimension of $\cM$.
\end{definition}

\subsection{Example} \label{example}
We give an example of calculation of the invariant defined 
in \S \ref{invariant}.
For preparation,
we show the following two lemmas.

\begin{lemma} \label{sum}
Let $M_i$ $(i=1,2)$ be an oriented closed $4$-manifold with $b^+(M_i)>1$ and
let $\cL_i$ be a spin$^c$ structure on $M_i$.
Assume that $(M_1,\cL_1)$ and $(M_2,\cL_2)$ satisfy the conditions $(*)$,
then $(M_1 \# M_2,\cL_1 \# \cL_2)$ also satisfies the condition $(*)$.
\end{lemma}

\begin{proof}
The condition $(*)_2$ is satisfied for $(M_1 \# M_2,\cL_1 \# \cL_2)$
by the definition of $c_{ij}$.
The condition $(*)_1$ is satisfied for $(M_1 \# M_2,\cL_1 \# \cL_2)$
by the sum formula of the index of the Dirac operator.
\end{proof}

We write $\Sigma_g$ for an oriented closed surface of genus $g$.

\begin{lemma} \label{K3}
Suppose  $M$ is a $K3$ surface
or $\Sigma_g \times \Sigma_{g'}$ with $g$ and $g'$ odd.
Let $\cL$ be a spin$^c$ structure on $M$ 
which is induced by the complex structure.
Then $(M,\cL)$ satisfies the condition $(*)$.
\end{lemma}

\begin{proof}
Note that $c_1(\det \cL)=-c_1(K_M)$.
Let $M$ be a $K3$ surface. 
The first Betti number of $M$ is equal to $0$, so the condition $(*)_2$ is satisfied.
By the index theorem \cite{AS},
the index of the Dirac operator is
\[
a=\frac{c_1(\det \cL)^2-\tau(M)}{8}
=\frac{0-(3-19)}{8}
=2
\equiv 0 \mod 2.
\]
Hence $(M,\cL)$ satisfies the condition $(*)$ when $M$ is a K3 surface.\\
Let $M$ be $\Sigma_g \times \Sigma_{g'}$ with $g$ and $g'$ odd.
Then we have
\[
c_1(\det \cL)=-c_1(K_M)=2(1-g)\alpha+2(1-g')\alpha'
\]
where $\alpha$ and $\alpha'$ are the standard generators of
$H^2(\Sigma_g;\Z)$ and $H^2(\Sigma_{g'};\Z)$.
Since $g$ and $g'$ are odd,
we have $c_1(\det \cL) \equiv 0 \mod 4$,
and then
\[
c_{ij}=\frac{1}{2}\< c_1(\det \cL) \alpha_i \alpha_j,[M] \>
\equiv 0 \mod 2,
\]
which implies the condition $(*)_2$.

By the index theorem,
the index of the Dirac operator is given by
\[
a=\frac{c_1(\det \cL)^2-\tau (M)}{8}
=\frac{c_1(\det \cL)^2}{8}.
\]
Because $c_1(\det \cL)^2 \equiv 0 \mod 16$,
we have $a \equiv 0 \mod 2$.
Hence the condition $(*)_1$ is satisfied.
\end{proof}

Let $M_j$ be a $K3$ surface or $\Sigma_g \times \Sigma_{g'}$,
where $g,g'$ are odd.
By Lemma \ref{sum} and Lemma \ref{K3},
the pair $(\#_j^l M_j,\#_j^l \cL_j)$
satisfies the conditions $(*)$,
where $\cL_j$ is a spin$^c$ structure on $M_j$ 
induced by the complex structure.
We show that the invariant $\sigma_{\#_{j=1}^l M_j}(\#_{j=1}^l \cL_j,\cO,L)$
is non-trivial when $l$ is $2$ or $3$.

\begin{theorem} \label{main thm}
Let $M_j$ be a $K3$ surface or $\Sigma_g \times \Sigma_{g'}$
with $g,g'$ odd
and $\cL_j$ be a spin$^c$ structure on $M_j$
which is induced by the complex structure.
Put $M=\#_{j=1}^l M_j$ and $\cL=\#_{j=1}^l \cL_j$
for $l=2$ or $l=3$.
Let $\sigma_M^0(\cL,\cO,L)$ be the image of $\sigma_M(\cL,\cO,L)$
under the natural map 
$\Omega_{l-1}^{spin}(T) \rightarrow \Omega_{l-1}^{spin}(*)$.
Then $\sigma_{M}^0(\cL,\cO,L)$ is non-trivial
in $\Omega^{spin}_{l-1}(*) \cong \Z_2$.
\end{theorem}

\begin{proof}
Let $L \rightarrow T$ be a square root of $\det \Ind(D)$.
If $l=2$, the dimension of the moduli space is one,
so the invariant $\sigma^0_{M}(\cL,\cO,L)$ is in 
the one dimensional spin bordism group $\Omega_1^{spin}(*) \cong \Z_2$,
and if $l=3$,
the invariant $\sigma^0_M(\cL,\cO,L)$ is in 
the two dimensional spin bordism group $\Omega_2^{spin}(*) \cong \Z_2$.
We will calculate the invariant for $l=2$ for simplicity.

Let $f_j:V_j \rightarrow W_j$ be
a finite dimensional approximation of the Seiberg-Witten map on $M_j$ 
such that $m_j=\dim W_{j,\C}$ is even,
and set $f=f_1 \times f_2:V=V_1 \times V_2 \rightarrow W=W_1 \times W_2$.
We make use of Bauer's construction (Theorem 1.1 in \cite{B}).
Bauer proved that there is a finite dimensional approximation on $M$
which is $U(1)$-equivariantly homotopic to $f$.

In general,
for a K\"ahler surface $M$ with $b^+(M)>1$ and 
a spin$^c$ structure $\cL$ on $M$ induced by the complex structure, 
the Seiberg-Witten moduli space $\cM_M(\cL,g,\eta)$ 
consists of smooth one point,
where $g$ is the K\"ahler metric and $\eta$ is a suitable $2$-form.
See, for example, \cite{N}.
Thus we may assume that
$\cM_j=f_j^{-1}(0)/U(1)$ is one point.
Hence $f_j^{-1}(0) \cong S^1$ and 
$\cM=f^{-1}(0)/U(1)_d=(f_1 \times f_2)^{-1}(0)/U(1)_d \cong S^1$,
where $U(1)_d$ is the diagonal of $U(1) \times U(1)$.
For some $t_j \in T_j=H^1(M_j;\R)/H^1(M_j;\Z)$,
$f_j^{-1}(0)$ lies in a fiber $V_{j,t_j}$ of $V_j \rightarrow T_j$.
Take a small open neighborhood of $t_j$ such that 
$V_j|_{U_j} \cong U_j \times \C^{m_j+a_j} \times \R^n_j$,
where $a_j$ is the index of the Dirac operator 
associated with $\cL_j$.
Set $S_j=U_j \times (\C^{m_j+a_j} \backslash \{ 0 \}) \times \R^{n_j}$
and $S=\prod_{j=1}^2 S_j$,
then $S$ has  a $U(1)_d$-action and a $U(1) \times U(1)$-action.
The $U(1)_d$-action is defined by the scalar product on
$\prod_{j=1}^2 (\C^{m_j+a_j} \backslash \{ 0 \})$.
And for $(\alpha_1,\alpha_2) \in U(1) \times U(1)$,
we define the action of $(\alpha_1,\alpha_2)$ on $S$ by
the scalar product of $\alpha_1$ on $(\C^{m_1+a_1} \backslash \{ 0 \})$ and
the scalar product of $\alpha_2$ on $(\C^{m_2+a_2} \backslash \{ 0 \})$.
Set $\bar{S}=S/U(1)_d$.

We write $\xi$ for a spin structure on $\bar{V}=V_{irr}/U(1)_d$ 
induced by $L$.
The restriction $\xi|_{\cM}$ of $\xi$ to $\cM$ is 
equal to $(\xi|_{\bar{S}})|_{\cM}$.
Since $H^1(\bar{S};\Z_2)=0$,
$\bar{S}$ has just one spin structure.
So it is sufficient to consider 
the restriction of the unique spin structure on $\bar{S}$ to $\cM$.

Put $U(1)_q=U(1) \times U(1)/U(1)_d \cong U(1)$,
then the $U(1) \times U(1)$-action on $S$ induces
a free $U(1)_q$-action on $\bar{S}$ and
$\bar{S}/U(1)_q=\bar{S}_1 \times \bar{S}_2$,
where 
$\bar{S}_j=S_j/U(1)
\cong U_j \times \CP^{m_j+a_j-1} \times \R_{>0} \times \R^{n_j}$.
Moreover this $U(1)_q$-action preserves $\cM \subset \bar{S}$
and induces a free $U(1)_q$-action on $\cM \cong S^1$.
Since $m_j+a_j-1$ is odd, $T\bar{S}_j$ has a spin structure.
So $T(\bar{S}/U(1)_q)$ has a spin structure.
Take a spin structure $\eta$ on $T(\bar{S}/U(1)_q) \oplus \underline{\R}$.
Let $p:\bar{S} \rightarrow \bar{S}/U(1)_q$ be the projection.
Then there is a natural isomorphism
$T\bar{S} \cong p^*(T(\bar{S}/U(1)_q) \oplus \underline{\R})$.
So $p^*(\eta)$ is the unique spin structure $\xi$ on $T\bar{S}$.
Because $p$ is the projection $\bar{S} \rightarrow \bar{S}/U(1)_q$,
the $U(1)_q$-action on $\bar{S}$ lifts to an action on $\xi=p^*(\eta)$.
So the $U(1)_q$-action on $\cM \cong S^1$ lifts to an action on 
restriction of $\xi|_{\cM}$.
In the same way,
we can prove that 
the $U(1)_q$-action on $\cM$ lifts to
an action on the spin structure on $\bar{E}|_{\cM}$.
Since
$f|_{S}=f_1|_{S_1} \times f_2|_{S_2}
:S_1 \times S_2 \rightarrow W_1 \times W_2$
is $U(1) \times U(1)$-equivariant,
the $U(1)_q$-action on $\cM$ lifts to
an action on the spin structure of $\cN$ induced by
$f$ and the spin structure on $\bar{E}|_{\cM}$.
Therefore the $U(1)_q$-action on $\cM$ lifts to
an action on the spin structure on $\cM$ induced by $f,\cO$ and $L$.
Such a spin structure determines
a non-trivial class in $\Omega^{spin}_1(*) \cong \Z_2$,
so $\sigma^0_{M}(\cL,\cO,L)$ is non-trivial class in $\Omega^{spin}_1(*)$
(See \cite{K}).

In the case of $l=3$, 
$\cM$ is the $2$-dimensional torus if we perturb the equations suitably.
We can show that the spin structure on $\cM$ is the Lie group spin structure as in  the case of $l=2$ and the spin bordism class $\sigma^0_{M}(\cL,\cO,L)$ is non-trivial in $\Omega^{spin}_2(*) \cong \Z_2$.
\end{proof}

\begin{remark}
Let $l$ be larger or equal to $4$.
Then we may assume that the moduli space is 
a $(l-1)$-dimensional torus $T^{l-1}$.
In the same way as in Theorem \ref{main thm},
we can see that 
the spin structure on $\cM$ induced by $f,\cO$ and $L$
is equal to the spin structure induced by the Lie group structure of $T^{l-1}$.
Such a spin structure is trivial in $\Omega^{spin}_{l-1}(*)$
if $l$ is larger or equal to $4$.
Hence $\sigma^0_M(\cL,\cO,L)$ is trivial in $\Omega^{spin}_{l-1}(*)$
when $l$ is larger than or equal to $4$.
\end{remark}

By Theorem $\ref{main thm}$,
we obtain the adjunction inequality for $M$.
See \cite{KM} for proof.

\begin{corollary}
Let $M_j$,$M$ and $\cL$ be as in Theorem \ref{main thm}.
Assume that an oriented, closed surface $\Sigma$ of positive genus is 
embedded in $M$ and
its self intersection number 
$\Sigma \cdot \Sigma$ is nonnegative. Then
\[
\Sigma \cdot \Sigma
\leq
\< c_1(\det \cL),[\Sigma] \>+
2g(\Sigma)-2,
\]
where $g(\Sigma)$ is the genus of $\Sigma$.
\end{corollary}

There are applications of Theorem $\ref{main thm}$
to computation of the Yamabe invariant and nonexistence of Einstein metric.

\begin{definition}
Let $M$ be an oriented, closed $4$-manifold.
Then the Yamabe invariant of $M$ is defined by
\[
\cY(M)=
\sup_{\gamma \in \Conf(M)}\inf_{g \in \gamma}
\frac{\int_M s_g d\mu_g}
{\left( \int_M d \mu_g \right)^{\frac{1}{2}}}
\]
where $\Conf(M)$ is the space of conformal classes of Riemannian metrics
on $M$, $s_g$ is the scalar curvature and $d\mu_g$ is the volume form of $g$.
\end{definition}

\begin{theorem} \label{Yamabe}
Let $M_j$ and $M$ be as in Theorem \ref{main thm},
and $N_1$ an oriented, closed, negative definite $4$-manifold admitting a Riemannian metric with scalar curvature nonnegative at each point.
Then
\[
\cY(M \# N_1)=-4\pi \sqrt{2\sum_{j=1}^l c_1(M_j)^2} .
\]
\end{theorem}

\begin{theorem} \label{non existence}
Let $M_j$ and $M$ be as in Theorem \ref{main thm}.
If $N_2$ be an oriented, closed, negative definite $4$-manifold satisfying
\begin{equation} \label{Ishida-LeBrun}
4l-(2\chi(N_2)+3\tau(N_2))
\geq
\frac{1}{3}\sum_{j=1}^l c_1(M_j)^2,
\end{equation}
then $M \# N_2$ does not admit an Einstein metric.
\end{theorem}

{\it Proof of Theorem \ref{Yamabe} and Theorem \ref{non existence}.}
In \cite{IL}, 
Ishida and LeBrun showed a similar statement
under a somewhat different assumption (Theorem D).
The main point of their argument is non-vanishing of
the Bauer-Furuta invariant.
In our case,
the invariant $\sigma_M(\cL,\cO,L)$ is non-trivial.
Hence we can apply their argument to our situation.
\qed \\

On the other hand,
there is a topological obstruction
for $4$-manifolds to have an Einstein metric (\cite{H}).

\begin{theorem}[Hitchin-Thorpe inequality \cite{H}]
Let $X$ be an oriented closed $4$-manifold admitting an Einstein metric, 
then
\begin{equation} \label{Hitchin-Thorpe}
3 |\tau (X) | \leq 2 \chi(X).
\end{equation}
\end{theorem}

\begin{example}
Let $M_i=\Sigma_{g_i} \times \Sigma_{g_i'}$ 
for positive odd integers $g_i,g_i'$,
let $M=M_1 \# M_2$ and let
$N=( \#^r \overline{\CP^2} ) 
\# ( \#^s S^1 \times S^3 )$.
Then $b^+(N)=0$ and 
the inequality (\ref{Ishida-LeBrun}) is satisfied
if $r \geq \displaystyle \frac{8}{3}G-4s-4$,
where $G:=\displaystyle \sum_{i=1}^2 (g_i-1)(g_i'-1)$.
By Theorem \ref{non existence},
$X=M \# N$ does not admit an Einstein metric
when $r \geq \displaystyle \frac{8}{3}G-4s-4$.
On the other hand,
if $r \leq 8G-4s-4$,
then $X$ satisfies the Hitchin-Thorpe inequality (\ref{Hitchin-Thorpe}).
Thus if
\[
\frac{8}{3}G-4s-4 \leq r
\leq 8G-4s-4,
\]
$X$ satisfies the Hitchin-Thorpe inequality,
but does not admit an Einstein metric.
\end{example}

\end{document}